\documentclass[9pt,twoside,a4paper]{article}
\usepackage{amsmath}
\usepackage{amssymb}
\usepackage{wrapfig}
\usepackage{graphicx}
\newcommand{\txtm}[1]{\mbox{$\smash{#1}$}}
\title{On mathematics and knowledge}
\author{{Daniel Canarutto} \\
{\small Department of Mathematics and Informatics, Florence University} \\
{\small ORCID iD: 0000-0002-8512-0014}}
\begin{document}
\bibliographystyle{alpha}
\maketitle 
\begin{abstract}
A sketch of some of the fundamental notions related to the nature of knowledge is offered,
with special focus on the role of mathematics and my own opinions.
No single idea exposed here is entirely original;
indeed, this topic has been explored by legions of philosophers,
mathematicians, and scientists throughout history|%
I listed a few related books among the references%
~\cite{Lorenz1973,Osinga1973,Gregory1981,Hofstadter1999,Sacks2008,Penrose2016,Hossenfelder2020}.
I look forward to observations and criticisms that may serve to improve this exposition.
\end{abstract}

\tableofcontents~\thispagestyle{empty}
\vfill\newpage

\section{On knowledge}\label{On knowledge}

\begin{flushright}
\emph{Considerate la vostra semenza:}\\
\emph{fatti non foste a viver come bruti,}\\
\emph{ma per seguir virtute e canoscenza.}\footnote{
Consider your lineage:
you weren't made to live like brutes,
but to follow virtue and knowledge.}
\smallbreak
Dante Alighieri, Inferno, Canto XXVI, 117--120.
\end{flushright}

First, I should make clear that I am not a fundamentalist of rational thought:
I think that there exist legitimate, meaningful questions
which cannot be answered rationally.
In fact, later I will also discuss rationality itself, and argue about its limits.

Nevertheless, the scope of this brief essay lies within the boundaries of science.
Accordingly, a good starting point is provided by the ideas of Konrad Lorenz~\cite{Lorenz1973}.

Human knowledge is based on an interactive process,
in which the knowing subject is confronted with the data of the surrounding world,
which are the object of knowing.
Actually, all living organisms are capable of knowing,
both by genetic inheritance and by learning.
Generally speaking,
we may say that this capability has naturally evolved during phylogenesis,
as it is essential for survival.
Indeed, any organism must `know' how to react in many different situations.
At the lowest level, this knowledge consists of a limited spectrum of innate reactions,
such as moving away from obstacles or selecting what to eat.
It is actual knowledge, however, suitable for thriving in the organism's habitat.

With evolution, the forms of knowledge have grown increasingly complex,
eventually generating the capability of learning from experiences,
and of conveying the acquired knowledge by non-hereditary means
(e.g.~by example, through language).
The advantages of such refinements are evident.\footnote{
Several years ago I made the following ethological observation
at Everglades National Park (Florida).
Here and there, on the roadside, there were remains of small animals crushed by the cars,
which attracted various birds such as herons, vultures, and crows.
When a car approached, herons and vultures would flee up the road in two or three short flights,
and finally get away.
The crows would just hop sideways and then return to their food as soon as the car had passed.}

We can observe in us, human beings, the whole evolutive spectrum of knowledge,
from the instinctive reactions of the newly-born to the performance of the most advanced minds.
Some of the latter are truly amazing.
Johann Sebastian Bach, for example,
was able to improvise a fugue for five voices on the keyboard,
based on a theme assigned to him on the spot.
Michelangelo imagined sculpted figures before `extracting' them from the stone.
There are chess players capable of playing dozens of games simultaneously blindly
(so I am told).

But the above examples are extreme cases of an ability that we all possess, to some extent,
and which we may call \emph{imagination}.
Inside our mind we can form, with variable accuracy and fidelity,
a virtual image of the world outside, or even of a non-existing world.
Within this image we can try and anticipate the future,
remember what is there in a certain location, and, most important,
perform mental experiments.
Just to give a trivial example, if a soccer ball got stuck in the branches of a tree,
I can imagine throwing a stone to make it fall,
and if I'm good enough then I may actually succeed.
Even some animals are known to be capable of doing such things, up to a point.

In order to be effective, imagination has to be trained through experience.
Indeed, the unending interaction and mutual feedback between experience and imagination is,
perhaps, the most powerful knowledge increasing mechanism.\footnote{
See also ne notion of `OODA loop' ({\tt https://en.wikipedia.org/wiki/OODA\_loop})
in the thought of John Boyd, examined in detail by F.P.B.~Osinga~\cite{Osinga1973}.}

A deeper understanding of imagination,
of its modalities and underlying physiological processes,
is a formidable task and certainly out of the scope of this paper,
but we can make a few basic observations.
The mind is equipped with the ability to discover regularities
in this mirror of the world that it possesses inside it;
this is obviously very important for coping with the environment
and increasing the probability of survival.
Furthermore, the mind can filter such regularities out of the different states of itself,
thus making it natural to assume that there is an exterior reality.\footnote{
Lorenz calls this inclination `hypothetic realism'.}
The discovery of regularities forms the basis for abstract thinking and mathematical thinking.

Indeed, mathematics can be essentially viewed as the \emph{description of regularities}.
Mathematical models are simplified representations of phenomena,
which greatly enhance our predictive ability.
A crucial point to remember, and that we shall further discuss in the sequel,
is that a mathematical model, or a theory,\footnote{
We are not making a sharp distinction between `model' and `theory';
roughly speaking, a theory can be thought of as more general and with a larger scope.}
always has a limited scope,
which can be explored and possibly expanded;
but the model and the reality should always be kept clearly distinct,
in order to avoid errors and missteps.
When a model fits especially well the phenomena under study then we tend to think:
``Aha! Now I understand!''

\section{On mathematics}\label{On mathematics}

Summarizing, Nature appears to be not a totally chaotic mess: it exhibits regularities;
mathematics is, basically, the description of regularities.
We must now make a distinction between usual mathematical thinking and rigorous formalisms.

In the most essential acception, a mathematical environment (a model, a theory...)
is a \emph{formal system},
consisting of a set of symbols (an \emph{alphabet}) and of a set of \emph{transformation rules}.
These have a purely operative meaning:
we do not define notions such as `set' or `rule';
we say that these notions are `primitive'.

Sequences of symbols are called assertions.
The rules specify how a given assertion
can be transformed into other assertions (the `consequences' of the former).
Such transformations are also called `deductions'.
Starting from given assertions (`axioms'),
any further assertion that can be deducted from them is a `theorem'.

A few observations are in order.
First, formal logic can be included into the considered formal system.\footnote{
A readable account of these matters can be found in Hofstadter~\cite{Hofstadter1999}.}
Second, when we fix a formal system and its axioms,
we have no guarantee that it is self-consistent.
Third, a sufficiently complex formal system is \emph{incomplete}:
there are \emph{undecidable} assertions, that is,
assertions that cannot be demonstrated to be either true or false
(neither the assertion nor its opposite are theorems).
Does this mean that the edifice of knowledge is flawed?
Or that we can know nothing, actually?
Not at all: remember,
we only expect that our formal model fits reality within a limited scope.
Experience allows us to trust the model as long as we stay within certain proven limits.
Of course, we will be exploring its mathematical properties and testing the limits,
but that's a separate activity from application.

Computer programs for symbolic computation 
can provide usable realizations of formal systems.
Usual mathematical thinking, however,
is somewhat different under several respects.\footnote{
This distinction may fade in the future as artificial intelligence progresses.}
One important observation is that 
there are no preferred derivation paths in a formal system,
and a mathematical expression has no preferred or `final' form.
A mathematician, instead, selects a path using ordinary language,
intuition, imagination, personal taste;
a mathematician has a purpose in mind,
be it build a better model or a neater theory.
Even a theorem's proof, which is supposed to be precise and rigorous,
is not presented as a formal derivation in the strict sense sketched above.
Indeed, managing a true formal derivation can be difficult or even impossible
for a human mind;
however, one expects that the proof, if it is correct,
can be cast in such a way at least in principle.
Euclid's geometry can be thought of as the prototype of a model 
that strives to adhere to the idea of a formal system.

Whatever the level of formality adopted,
a mathematician's choices also rely on intuitive meanings assigned to mathematical concepts;
the symbols by themselves have no meanings.
For example, consider the script
\[ \begin{cases}a\in A\\ f:A\to B\end{cases}\quad\Rightarrow\quad f(a)\in B~. \]
In the abstract context, \txtm{a\in A} and \txtm{f:A\to B} are purely formal assertions
(here assumed jointly), as well as the ensuing theorem \txtm{f(a)\in B}\,.
On the other hand, saying `$a$ is an element of the set $A$' is an interpretation,
which is useful for the the human mind to handle the symbols and rules through visualization,
and to select a deductive path.
This selection process, not deterministic, may be called `induction'.

Similarly, consider the scripts
\begin{align*}
& \lim_{n\to\infty}a_n = a~,
\\
&\forall\epsilon>0~\exists n_\epsilon\in\mathbb{N}\,:\, 
n>n_\epsilon\Rightarrow|a-a_n|<\epsilon~.
\end{align*}
The left-hand side of the former script should be actually seen as a unique symbol,
but contains glyphs that are read as `limit', `infinity', `tends~to';
the latter script is an assertion that contains symbols
that are read `for~all', `exists', `such~that', `follows', and so on.
Indeed, this assertion can be read in human language and, as such,
is somewhat easier to grasp.
But `exists' or `infinity' are basically mere symbols, to be treated according to certain rules.
The intuitive meanings we assign to them guide our steps in the formal world.\footnote{
My daughter, when she was about three years old, once asked me:
``Dad, do the numbers never end?''}

The distinction between formal rules and intuitive notions is especially relevant
in practical applications of mathematics.
Clarity in physical theories, in particular,
requires to distinguish the abstract mathematical model,
the description of experiments,
and the links between these two sides.
Thus, at all levels, a theory is basically characterized by \emph{operative definitions}.
Note how this operativity aspect is also true of the formal mathematical model
and of mathematics itself.

\section{The beauty of mathematics}
\label{The beauty of mathematics}

As we try and build a physical theory, research proceeds in two different directions,
that we might label `external' and `internal'.
Externally, we modify our model's assumptions in order to fit the phenomena better and better;
internally, we extend and refine the mathematics.
The latter process may become, in part, autonomous from the practical applications.

Let us try and understand what drives research in pure mathematics.\footnote{
Apart from the obvious 'career hunger' that drives all research.}
After all, in that field there is no limit or rule as to what is worth pursuing.
However the history of mathematics seem to follow a path of \emph{discoveries},
along which certain milestones stand out as beautiful truths.
I'm thinking of the constructions of the real and complex numbers,
of analysis, differential geometry, and so on.
When we learned about all this, as students, we were gratified by aesthetic appreciation.
Perhaps, this kind of appreciation is not that far
from what one feels for artistic masterpieces.
Take, for example, the great works of Bach, Mozart, Beethoven, and many others;
there is substantial agreement about their value,
at least among those who do appreciate music.
A similar observation is valid for other forms of art.

This implies, in my opinion, that there is \emph{truth} in our aesthetic judgments.
And now the question: \emph{where} do truth and beauty reside?
It seems to me that there are essentially \emph{two} possible answers.
One answer is that, in some misterious way, they reflect the inner workings of our minds.
The neurologist and popularizer O.~Sacks observed~\cite{Sacks2008} how
``for virtually all of us, music has great power, whether or not we seek it out
or think of ourselves as particularly ``musical''.
This propensity to music|this ``musicophilia''|shows itself in infancy,
is manifest and central in every culture,
and probably goes back to the very beginnings of our species.
...
While birdsong has obvious adaptative uses...
The origin of human music is less easy to understand.
Darwin himself was evidently puzzled...''.
Sacks goes on by quoting various proposals that have been advanced
in order to frame music in the Darwinian context,
noting that none is definitively convincing.
Moreover, the clinical cases presented in the book tend to show how music
is within us at a very deep level.

The other possible answer to the above question is that truth and beauty are absolute qualities,
which came to us from some place outside, like Plato's Hyperuranion.
This point of view may not seem very `rational', or `scientific',
but has been seriously discussed by mathematicians and scientists at the highest levels.
Alexander Grothendieck, in particular,
wrote that mathematics belongs to the nature of God.\footnote{
``Il convient cependant de faire exception ici des lois math\'ematiques.
Ces lois peuvent \^etre d\'ecouvertes par l'homme, mais elles ne sont cr\'e\'es ni par l'homme,
ni m\^eme par Dieu. ... Je sens les lois math\'ematiques comme faisant partie
de la nature m\^eme de Dieu ...''
(Grothendieck~\cite{Grothendieck}, p.\ 100 ch.\ 31).}

As for Paul Dirac, his claim  that
``it is more important to have beauty in one's equations than to have them fit experiment''
is well-known~\cite{Morrison2016}.
It is interesting to note that, though he clearly expressed atheistic feelings,
he also said: ``It seems to be one of the fundamental features of nature
that fundamental physical laws are described in terms of a mathematical theory
of great beauty and power, needing quite a high standard of mathematics for one to understand it.
You may wonder: Why is nature constructed along these lines?
One can only answer that our present knowledge seems to show that nature is so constructed.
We simply have to accept it.
One could perhaps describe the situation by saying that God
is a mathematician of a very high order,
and He used very advanced mathematics in constructing the universe.
Our feeble attempts at mathematics enable us to understand a bit of the universe,
and as we proceed to develop higher and higher mathematics 
we can hope to understand the universe better.''

I would say that the sentences quoted above
place Dirac in the ranks of idealist thought, of Platonic legacy.
His putting beauty beyond fitting experiment may be justified by the observation
that a preliminary misfit may not be a good reason
for the immediate rejection of a fine theory.
But now let me ask, what is `beauty in one's equations'?
When a new, ground-breaking theory is formulated,
it may generate a strong impression,
related to subtle paradygm shifts and unexpected mathematical background,
but it is not just a matter of `equations'.
After all, if you set out to write down differential equations,
the possibilities are endless.
You need a guide and some selection criterion.
Consider some of the great physics revolutions of the twentieth century:
General Relativity, Quantum Theory, and Dirac's theory of spin.
Their mathematical backgrounds are, basically, curved spacetime, Hilbert spaces, and spinor spaces;
these are all \emph{geometric} notions,
although not pertaining to the usual geometry of our perception.
When we select a geometrical background,
the natural equations become essentially fixed.
On the other hand, geometry is deep inside our minds,
as it does play an an essential role in the survival kit of our species.\footnote{
As evidenced, in particular, by the structure of the inner ear dedicated to spatial orientation.}
Moreover, strange or unusual geometries tend to exert a great fascination on speculative minds.

Esthetic appreciation, and the curiosity for the discovery of new abstract relations,
are the bases for intuition and creativity;
eventually, they power concrete scientific progress.\footnote{
`The intuitive mind is a sacred gift and the rational mind is a faithful servant'
(attributed to A. Einstein).}
This aspect is essential
in explaining the evolutive success of the human species.
Hence it is pointless to ask ``what is this abstract research for?''
If we already knew, there would be no discovery.

Eventually, however, aesthetics should not be given the last word:
the validity of aesthetic judgments is different from the truth of a physical theory.
There are pitfalls hiding behind the idealist way of thinking,\footnote{ 
In particular, see Hossenfelder~\cite{Hossenfelder2020}.}
and I will discuss some of them in the next sections.

\section{Empiricism and idealism}\label{Empiricism and idealism}

\begin{flushright}
\emph{When you have eliminated all which is impossible,}\\
\emph{then whatever remains, however improbable, must be the truth.}
\smallbreak
Sherlock Holmes
\end{flushright}

The ancient dispute between idealism and empiricism
runs through the entire history of philosophy in various forms.
Some of the early premises of empiricism are untenable today;
in particular, it is impossible for us to see the mind as a \emph{tabula rasa} (a blank slate),
that contains nothing to begin with and acquires knowledge only through experience.
Nevertheless, let us keep the term `empiricism' as signifying the opposition to the view
that true scientific discoveries can be made by working only inside our minds.
In this limited sense, I would call myself an empiricist.\footnote{
The empirical method is popularized by investigators in detective stories.
Inspector Maigret uses to linger around,
listening to people and trying to absorb the atmosphere of the places,
until a possible unfolding of the events forms inside his mind.
Nero Wolfe never leaves his house, preferably,
and lets Archie Goodwin bring in the facts;
eventually he devises a model that fits them.
As for Sherlock Holmes' famous sentence,
it could be used as a manifesto of empiricism|%
though always finding the solution is an illusion.}

The observations offered in \S\ref{On mathematics} lead us to exclude that mathematics,
and the scientific research that is based on it,
can be seen (as in the popular view) as totally precise, rational activities|%
whatever it means.
Actually, we can never be certain about the internal consistency of any mathematical model.
We just try to apply the considered model in a specific context.
We gain trust in the model's validity, and probe its limits,
by comparing its predictions with experiments.

Let me say it again: the purely abstract work is an essential part of scientific research;
furthermore, when a mathematical model fits the experiments,
looking for extensions and generalizations, in order to expand its scope,
is certainly natural and opportune.
Indeed, this `inductive' process can produce breakthroughs;
but it can also lead to dead ends.
The pitfall consists of presuming that such generalitions say something about Nature
by their own virtue.
This can be the case, sometimes, but it is not at all guaranteed.\footnote{
No doubt, the reader will be able to think of examples of situations
where generalization has indeed worked
(just to name one, the extension of Lagrangian theory
from the motion of a system with finitely many
degrees of freedom to classical and quantum field theory).}
Furthermore, some of the idealist approaches of today seem to dispense
with a criterion that is, in my opinion, of the utmost importance, that is
Occam's razor principle: \emph{entities are not to be multiplied beyond necessity}.
Granted that abstract exploration is needed,
a selection must be eventually made, based on experience;
otherwise we could drift indefinitely.

\smallbreak
Next I will discuss some cases, from ancient Greece to present times.

\smallbreak\noindent
$\bullet$\quad To start right from the beginning,
consider Zeno's argument for demonstrating the impossibility of motion:
for the arrow to hit the target, it must first travel half of the way,
then half of the remaining half, and so on;
as it has to travel an infinity of intervals, it will never arrive.
The solution to this `paradox', as it is offered to highschool students,
uses  the observation that \hbox{$\frac12+\frac14+\frac18+\dots=1$},
so that the arrow actually travels the whole distance
and does so in a finite time|%
since each interval is traveled in a time proportional to its length.
Now, certainly Nature does not perform such calculations (or `renormalizations').
The actual flaw in Zeno's argument lies in presuming that a purely mental construction
can be resolutive and reliable.

\smallbreak\noindent
$\bullet$\quad One of the clearest examples of an idealistic attitude
is the fascination exercised for centuries by the five regular polyhedra
or Platonic solids, that Plato himself associated with the elements.
We no longer regard them as candidates for explaining the structure of the universe,
but certain suggestions at fundamental theories may seem surprisingly similar
as for basic attitude.
Similarly, the geocentric description of the solar system before Kepler was based
on the circle, as the most perfect geometric figure.
In order to do astronomical calculations they had to consider complicated movements of epicycles,
but, eventually, the calculations were quite precise.
This observation could be used to contrast the opinion that quantum field theory
is all right as it is, since eventually, with the aid of renormalization,
yields precise calculations.

\smallbreak\noindent
$\bullet$\quad A delicate issue regards the role of the `real' numbers,
which are actually based on a highly sophisticated and abstract construction.
Just think about this:
the bulk of the field $\mathbb{R}$ of reals consists of elements
that cannot be expressed or characterized in any way by finite information;
the real numbers we deal with in practice constitute a very small subset of measure zero;
all the rationals are contained in an open subset of arbitrary small measure.
At the same time, there exist subsets, such as the Cantor set,
that are not countable but have measure zero.

Though all this defies our intuitive capacities,
the reals `exist' in the mathematical sense
which we hinted at in \S\ref{On mathematics}.
The reason why they turn out to be so important is that in $\mathbb{R}$
we can introduce a suitable notion of limit,
and prove all those theorems which make calculus consistent.
In other terms, the reals provide us with a powerful and sound context for calculations.
This is not the same as saying that they are to be included
into the fundamental notions of physics.
I would rather say that taking the reals for granted,
and including differentiable manifolds in the basic setting of a fundamental theory,
amounts to starting from strong, involved assumptions.

\smallbreak\noindent
$\bullet$\quad Newton's laws constitute one of the earliest big successes of theoretical physics.
That success, and especially the related deterministic paradigm
and the mirage of exact solutions,
has conditioned all the following development of science,
even outside of physics.
That conditioning was only partially removed with introduction of
quantum mechanics and, later, deterministic chaos;
later, I will argue that it is still present today and may be hindering the efforts
towards a new fundamental theory.

\smallbreak\noindent
$\bullet$\quad The case of probability theory is also interesting.
Probability can be, and is, also used at an elementary level.
In some expositions we may encounter a question that, in my opinion, is ill-posed:
what is the `definition' of probability?
Do we prefer the `subjective' or the `frequentist' approach?
The answer, of course, is that these are not definitions:
probability theory is a mathematical model, a part of the theory of measure,
that can be applied to different situations such as betting
or the study of the frequences of outcomes of experiments.

\smallbreak\noindent
$\bullet$\quad About quantum mechanics, a few observations come to mind.
Before its introduction, Niels Bohr looked forward to a `rational mechanics of the atom',
meaning a formal setting apt to describing microscopic phenomena;
and the foundations of the theory were actually laid down to that end.
Again and again, the validity of the theory has been confirmed by experiment
throughout a vast range of phenomena,
but some philosophical questions have arised.
Great attention has be devoted to the `interpretation' of quantum mechanics,
though the diverse interpretations that have been proposed do not appear
to be distinguishable through any empirical evidence.
Certain mathematical concepts of the formal setting seem to have acquired a status
of fundamental entities,
implying that any fundamental theory must be formulated in terms of these concepts.
This has had far-reaching consequences for theoretical physics,
which has drifted in immense oceans in the quest for a new theory
merging gravitation and quantum physics,
as the experimental evidence apt to falsify hypotheses became scarcer and scarcer.

\smallbreak\noindent
$\bullet$\quad One approach to fundamental physics I am not convinced of
is the approach called `quantization',
that strives to find a general method for deriving quantum physics from classical physics.
But there are several other problematic or unclear ideas floating around:
we may read about `the wave equation of the universe',
discuss whether Schr\"odinger's cat actually esists in a mixed state,
or whether quantum states are real or not.
Furthermore, by sticking to the original form of the Golden Rules of quantum mechanics,
the states actually observed take a somewhat secondary role,\footnote{
A clearer formulation can be achieved by replacing the Hilbert space
with a so-called \emph{rigged Hilbert space}~\cite{Canarutto2020a}.}
while by assigning a fundamental role to the algebra of observables
one is led to consider \emph{algebraic} approaches in which such common notions
as spacetime points should arise as derived notions.
I am not aware of any actual practical breakthrough in these theories,
nor in the great sea of string theories and beyond~\cite{Penrose2016}.
Though I am no expert,
I would bet we will see no experimentally verifiable progress of in any of these directions.
I await to be proved wrong...

\smallbreak\noindent
$\bullet$\quad Unification is one inductive procedure
which has been successfully employed more than once.
Maxwell's theory derives from the unification of the electric and magnetic fields.
It is no wonder that such accomplishment has been inspirational.
The following electroweak unification confirmed the validity of the unification procedure.
Further attempts were not as successful.
In particular, the idea of unifying the gravitational field with the matter fields
may remain a chimera.

\section{Particles, fields, and geometry}\label{Particles, fields, and geometry}

Let us return to the issue of the role of the real numbers,
which are associated with the notion of a continuum.
While the actual theory of the reals involves somewhat esoteric concepts,
the quest for the continuum stems from natural constructions such as dividing into equal parts
(the rational numbers) and the requirement of completeness (limits).
Furthermore, there is little doubt that the real numbers have been \emph{discovered}.
Hence, if we stand with Grothendieck in saying that mathematics is in the nature of God,
we might well assume that the Creation is based on the reals.

An argument I read once\footnote{Unfortunately I cannot find the precise reference any more.}
went more or less like this:
space and time can no more be considered a-priori categories of knowledge,
as intended by Kant;
hence, they must be replaced by a new category: the continuum.
If we start from this assumption, then it seems that field theories must be fundamental.
Frankly, I cannot agree with the view that the continuum must be regarded
as an a-priori category of knowledge.
Anyhow, from the empiricist point of view,
we may ask: is this assertion falsifiable?

John Baez~\cite{Baez2020} writes:
``...in every major theory of physics,
challenging mathematical questions arise from the assumption that spacetime is a continuum.
The continuum threatens us with infinities.
Do these infinities threaten our ability to extract predictions from these theories|%
or even our ability to formulate these theories in a precise way?
We can answer these questions, but only with hard work.
Is this a sign that we are somehow on the wrong track?
Is the continuum as we understand it only an approximation to some deeper model of spacetime?''

The refutation of the continuum hypothesis then requires a new theory,
based on discrete concepts,
in which the calculations of quantum field theory can be replaced
by more straightforward procedures,
somewhat similarly to the transition from the geocentric
to the eliocentric descriptions of the solar system.
Note how, still today, we are deeply entangled essentially with the same old issue
that Niels Bohr thought he settled with his complementarity principle:
field, or particle?

According to the continuum view, particles are an `epiphenomenon' and
the fundamental theory is a field theory.
Others think that the foundations of physics are discrete,
and that the continuous description is just a useful tool at some level.
The former approach is fairly standard in theoretical physics,
where one works out the `quantization' of classical fields.\footnote{
See for example the introduction of Weinberg's book~\cite{Weinberg1996}.} 
Now, the discrete aspects in the physics of quantum particles are evident,
so a theory founded on the continuum has the problem of explaining
how discreteness actually arises.
Tentative arguments draw analogies from established physics where discreteness arises,
e.g. the `vibration modes' of a bounded system (for example, a string),
or spectra of operators.

On the other hand, the discrete approach has his supporters.
A clear-cut version of this philosophy has been expounded by Penrose~\cite{Penrose1971}:
\begin{center}
\begin{minipage}{370pt}
``The idea is to concentrate only on things which \emph{are} discrete in
existing theory and try and use them as \emph{primary concepts}---then to build
up other things using these primary concepts as the basic building blocks.
Continuous concepts could emerge in a limit,
when we take more and more complicate systems.\\
...\\[6pt]
The central idea is that \emph{the system defines the geometry}...
The notion of space comes out as a \emph{convenience} at the end.''
\end{minipage}
\end{center}
Various proposals for a discrete approach
have also appeared more recently in the literature~\cite{Sorkin2003,Verlinde2010,Rovelli2010}.

\smallbreak
Basically, I regard the dichotomy between the continuum hypothesis
and the discrete hypothesis, or between fields and particles as fundamental concepts,
as anchored to the dichotomy between idealism and empiricism.
I deem I presented the reasons for this interpretation sufficiently clearly.
I will only add that my own choice is decidedly for empiricism and discreteness;
in a previous paper~\cite{Canarutto2014},
I presented a more detailed view of my ideas on this matter.

\newpage
\section{Conclusions}\label{Conclusions}

\begin{flushright}
\emph{There are more things in heaven and earth, Horatio,}\\
\emph{Than are dreamt of in your philosophy.}
\smallbreak
William Shakspeare, Hamlet
\end{flushright}

There is no precise scientific definition of `scientific method'.
There is no rational definition of rationality.
Nevertheless, making the distinction between scientific and unscientific
is of utter importance.

The scientific method can be loosely defined
as `an empirical method for acquiring knowledge',
characterized as a never-ending cycle including steps such as
observation, hypotheses, experiments.\footnote{
{\tt{https://en.wikipedia.org/wiki/Scientific\_method}}\,.}
In \S\ref{On knowledge} we observed that a feedback loop
is at the basis of essentially any knowledge-acquiring mechanism.
Let me remark that it is not just a simple loop in general,
as the mechanism may involve multiple sub-loops,
internal elaboration, and continuous comparison with the environment's responses.

Furthermore, this process is highly individual.
It depends on play, creativity, culture, previous knowledge.
Indeed, `understanding' \emph{is} subjective.
One scientist may be content with being able to make precise, trustable computations,
another may look for a rigourous mathematical model;
yet another my disagree with the former's conception of rigour.

Nevertheless, we may ask: when is a hypothesis scientifically confirmed?
Karl Popper's answer is univocal: never, the sole purpose of scientific research
is to \emph{falsify} hypotheses.
Adopting a less extreme position,
we might say that, according to our empiricist point of view,
a theory is not required to represent the ultimate foundations of Nature;
we are content with using it as a trustable description of Nature's workings
within a given scope.
Whether that scope is expandable or not, that remains to be seen.

Today we are often confronted with unproved hoaxes,
and scientists seem, sometimes, unable to contrast them.
But what gives us the guarantee that a scientist's opinion
is scientific and trustable?
In principle, nothing.
A pseudo-scientist might fill a paper of mathematical formulas
or other seemingly scientific reasoning,
and only produce nonsense.

But then, what is scientific reasoning?
A scientific answer to this question would be a logical paradox,
but this does not imply that we cannot distinguish.
Being able to distinguish and select between different answers
is a fundamental capacity of the mind,
related to the qualities that are called `intuition' and `wisdom' in common language,
though not so precisely characterizable in general.
I would rather say that scientific thought can be described
as a psychological attitude related to honesty,
and includes the capacity of disregarding personal convenience and preferences,
as well as unsubstantiated beliefs.

On the other hand, there are `fundamentalists of rationality'
who claim to decide, in principle, which matters are scientific and which are not;
the ground, here, is quite treacherous.
One way to put such position consists of saying that the only questions that make sense
are those that can be answered rationally.
In my opinion, this statement is meaningless,
as it can provide no actual selection criterion.

Many atheistic intellectuals, however, recognize that there is a fundamental mystery
that it is legitimate to question: the very existence of the universe.
Some believe that even this can be answered in scientific and rational terms,
but in practice all attempts of this type turn out to be logically flawed,
thus ending up confirming the opposite point of view.\footnote{
Some of these attempts have to do with: 
the use of the `anthropic principle';
the `theories' on the birth of our universe from something pre-existing;
the ruminations on the `probabilities of the axioms of mathematics'.}
The fundamental mystery could be definitively unfathomable by reason.
But over the centuries refined minds have questioned themselves|%
on various issues and in different ways|expressing profound thoughts
that cannot be dismissed so easily.
It is not true that current thinking cancels or is definitively superior
to anything that has been expressed in the past.

\vfill\newpage


\end{document}